% sj107.tex
\input amstex   %version 2.1
\documentstyle{amsppt}  %version 2.1
%\documentstyle{sjmod} :
%\begingroup
% \catcode`\@=11
 % annan kursivering av referenser:
\let\AMSpaper\paper
\let\AMSjour\jour
\let\AMSinbook\inbook
\def\paper{\AMSpaper\rm}
\def\jour{\AMSjour\it}
\def\inbook{\AMSinbook\it}
%\endgroup
%
\def\today{\ifcase\month\or January\or February\or March\or April\or May\or
   June\or July\or August\or September\or October\or November\or December\fi
   \space \number\day, \number\year}
\begingroup
  \count255=\time
  \divide\count255 by 60
  \count1=\count255
  \multiply\count255 by -60
  \advance\count255 by \time
  \ifnum \count255 < 10 \xdef\klockan{\the\count1.0\the\count255}
  \else\xdef\klockan{\the\count1.\the\count255}\fi
\endgroup
%
% inga sidhuvuden
\leftheadtext{\headtext}
\rightheadtext{\headtext}
    %ensidigt
\define\headtext{\relax} % om ej \draft
\define\draft{
  \gdef\headtext{\jobname.tex\quad DRAFT \today, \klockan}
  \footline={\vtop{\vskip-24pt\eightrm\hfil\headtext\hfil}}}    %pga sid 1
%
%roster romersk

% \qed at end of line
\redefine\qed{\begingroup\parfillskip=0pt\hbox{}\hfill\quad$\square$\par
                                                                  \endgroup}
% sidstorlek mm
\magnification\magstephalf\message{\string\magstephalf}
\pagewidth{150truemm}\pageheight{240truemm}\hoffset=3truemm %A4
\TagsOnRight

\define\Remark{\example}
\define\endRemark{\endexample}
\def\SJ#1{\csname SJ\romannumeral#1\endcsname}
\newcount\refnr \refnr=0
\define\refdef#1{\advance\refnr by 1\xdef#1{\the\refnr}}

%\define\wickl{\,\mathord:}
%\define\wickr{\mathord:\,}

  \def\eqd{\,\overset \text{d} \to=\,}
  
  \def\Var{\operatorname{Var}}
  \def\E{\operatorname{E}}
  \def\P{\operatorname{P}}
  
  \def\Bi{\operatorname{Bi}}
  \def\Be{\operatorname{Be}}
  \def\Po{\operatorname{Po}}
  \define\Ge{\operatorname{Ge}}
  \define\Exp{\operatorname{Exp}}
\define\0{{\rm 0--1}}
\def\kappa{\varkappa}

%\draft

\topmatter
\title
Large deviation inequalities for sums of indicator variables
\endtitle
\author Svante Janson\endauthor
\date\vbox{%\centerline{2 September 2016}\smallskip
\noindent
\it
This paper was written in  1994, but was never published because
I had overlooked
some  existing papers containing some of the inequalities.
Because of some recent interest in one of the inequalities, which does not seem
to be published anywhere else, it has now been lightly edited
and made available here. 
2 September, 2016
}
\enddate
\address Department of Mathematics, Uppsala University, PO Box 480,
S-751 06 Uppsala, Sweden
\endaddress
\email svante.janson\@math.uu.se \endemail
\thanks Supported by the G\"oran Gustafsson Foundation for
Research in Natural Sciences and Medicine
\endthanks

\abstract
A survey is given of some Chernoff type bounds  for the tail probabilities
$\P(X-EX\ge a)$ and $\P (X-EX\le a)$ when $X$ is a random variable that can
be written as a 
sum of indicator variables that are either independent or negatively related.
Most bounds are previously known and some comparisons are made.
\endabstract

\endtopmatter
\document

\head
1.
Introduction and conclusions
\endhead
The purpose of this paper is to give a survey of some simple upper bounds
for the 
probabilities $\P (X-EX\le -a)$ and $\P (X-EX\ge a)$, where $X$ is a random
variable that can be written as a sum $I_1+\dots+ I_n$ of \0 (indicator)
random variables. We consider both independent  and dependent variables $I_i$
(with strong restrictions in the dependent case).
Many of the inequalities extend to sums of more general bounded variables,
but we consider for simplicty only the indicator case.

Most of the bounds are known, see in particular 
Bennett (1962) and Hoeffding (1963), 
but are included for comparison and (partial) completeness. A few
versions seem to be new.
Many of the inequalities appear in various places, 
for example Janson, \L uczak and Ruci\'nski (2000), Chapter 2.
See also the book Boucheron, Lugosi and Massart (2013) which presents several
of these bounds and many extensions to other situation.

\subhead
Independent identically distributed summands
\endsubhead
The simplest case is when the indicator variables $I_i$ are independent and
identically distributed, $I_i\sim \Be (p)$, with $0<p<1$ (avoiding trivial
cases); then $X$ has the binomial distribution $\Bi (n,p)$. This case has been
studied by many authors, giving bounds or asymptotic results (sometimes in
greater generality); see for example Khintchine (1929), Cram\'er (1938), Feller
(1943), Chernoff (1952), Bahadur and Rao (1960), 
Bennett (1962), Hoeffding (1963), Littlewood (1969),
and the further references given in these papers.

We are here interested in explicit bounds for finite $n$ rather than
asymptotic results. One simple but powerful such bound was given by Chernoff
(1952); since for every $t\ge0$,
$$
\P (X\ge \E X+a)\le e^{-t( \E X+a)} \E e^{tX} \tag1.1
$$
and $Ee^{tX}=(1+p(e^t-1))^n $, we obtain by simple calculus, letting
$\lambda=\E X $ and assuming $0\le a\le n-\lambda$,
$$
\align\hskip-2em
\P (X\ge \E X+a) &\le \inf_{t\ge0}\exp \bigl(-at-npt+n \ln (1+p(e^t-1))\bigr)\\
 &= \exp\Bigl(-(\lambda+a)\ln \frac{\lambda+a}{\lambda} -(n-\lambda-a)\ln
 \frac{n-\lambda -a}{n-\lambda}\Bigr).\hskip-2em\tag1.2
 \endalign
 $$
 Similarly, for $0\le a\le \lambda$,
 $$
 \align\hskip-2em
 \P (X\le  \E X-a)&\le \inf_{t\ge0}\exp \bigl(-at+npt+n\ln (1+p(e^{-t}-1))\bigr)\\
& =\exp \Bigl(-(\lambda-a)\ln \frac{\lambda-a}{\lambda}-(n-\lambda+a)\ln
 \frac{n-\lambda+a}{n-\lambda}\Bigr).
\hskip-2em\tag1.3
\endalign
$$
 Chernoff (1952) proved also that the estimates (1.2) and (1.3) are
 asymptotically sharp in the sense that if $n\to\infty$ with $a/n$ and $p$
 fixed, then equality holds within factors $1+o(1)$ in the exponent.

 Simpler but (slightly) cruder bounds are easily obtained by finding suitable
 upper bounds for the right hand sides of (1.2) and (1.3), see for example
 Alon and Spencer (1992), Appendix A. We state here some, more or less
 well-known, such bounds. Proofs are given in Section 2.

\proclaim{Theorem 1}
Suppose that $X\sim \Bi (n,p)$ and let $q=1-p$. Then, for every $a\ge0$,
$$
\align
\P (X\ge  \E X+a) &\le \exp \biggl(- \frac{a^2}{2(npq+a(q-p)/3)} \biggr)
\le \exp \biggl(-\frac{a^2}{2(npq+a/3)}\biggr)\\
& \le \exp \Bigl(-\frac{a^2}{2npq }\Bigl(1-\frac{a}{3npq}\Bigr)\Bigr),\tag1.4\\
\P (X\ge  \E X+a)&\le \exp \biggl(-np\Bigl(\Bigl(1+\frac{a}{np}\Bigr)\ln
\Bigl(1+\frac{a}{np}\Bigr)-\frac{a}{np}\Bigr)\biggr),\tag1.5\\
\P (X\ge  \E X+a) &\le \exp \biggl(-\frac{a^2}{2np(1+a/3np)}\biggr)
\le \exp \Bigl(-\frac{a^2}{2np}\Bigl(1-\frac{a}{3np}\Bigr)\Bigr), \tag1.6\\
\P (X\le  \E X-a) &\le \exp \biggl(- \frac{a^2}{2(npq-a(q-p)/3)}\biggr)
\le \exp \biggl(-\frac{a^2}{2(npq+a/3)}\biggr)\\
& \le \exp \Bigl(-\frac{a^2}{2npq }\Bigl(1-\frac{a}{3npq}\Bigr)\Bigr),\tag1.7\\
\P (X\le  \E X-a)&\le \exp \biggl(-np\Bigl(\Bigl(1-\frac{a}{np}\Bigr)\ln
\Bigl(1-\frac{a}{np}\Bigr)+\frac{a}{np}\Bigr)\biggr),\tag1.8\\
\P (X\le  \E X-a) &\le \exp \Bigl(-\frac{a^2}{2np}\Bigr).\tag1.9
\endalign
$$
Moreover, if\/ $0\le p\le 1/2$, then
$$
\P (X\le  \E X-a)\le \exp\Bigl(-\frac{a^2}{2npq}\Bigr).\tag1.10
$$
\endproclaim
\Remark{Remark 1}
The estimates (1.4) and (1.7) (as well as (1.2) and (1.3))
are obvious ``mirror images''
of each other, and are equivalent by the substitution $X \to n-X$.
On the other hand, the remaining estimates in Theorem 1 are asymmetric, and
are useful mainly when $p$ is small.
\endRemark
\Remark{Remark 2}
Note that (1.5), (1.6), (1.8) and (1.9)
use $n$ and $p$ only in the combination $np=\lambda=\E X$.
\endRemark
\Remark{Remark 3}
It may also be observed that if we replace $np$ by $\lambda$,
then (1.5), (1.6), (1.8) and (1.9) are (e.g\. by continuity)
valid also when $X\sim\Po(\lambda)$; in fact,
(1.5) and (1.8) then become the Chernoff bounds for the Poisson distribution.
These Chernoff bounds too are asymptotically sharp in the sense, considering
for example (1.5),
that if $c>0$
is fixed and $c_1>(1+c)\ln(1+c)-c$, then
$\P (X\ge  \E X+c\lambda)> \exp (-c_1\lambda)$
for $X\sim\Po(\lambda)$ with
$\lambda$ large, and thus also for some $X\sim\Bi(n,p)$ with $np$ large and $p$
small.
In particular this implies that the simple bound (1.9) is {\it not} valid for
$\P (X\ge  \E X+a)$. This implies further, by considering $n-X$, that
(1.10) cannot hold without some restriction on $p$.
\endRemark

\subhead
Independent summands with different distributions
\endsubhead
The Chernoff bounds given above for the binomial distribution are easily
extended to the case when the \0 variables $I_i$ are independent with
different distributions, $I_i\sim \Be (p_i)$. In fact, as is well-known
(see for example Alon and Spencer (1992), Appendix A), if $\lambda= \E X
=\sum_1^n p_i$,  $p=\lambda/n$ (the average of $p_1,\dots,p_n)$, and we let
$X_0\sim \Bi (n,p)$ be a binomially distributed random variable with the same
$n$ and expectation as $X$, then by Jensen's inequality for the convex
function $x\mapsto-\ln \bigl(1+x(e^t-1)\bigr),$
$$
 \E e^{tX}=\prod _1^n \bigl(1+p_i(e^t-1)\bigr)
 \le \bigl(1+p(e^t-1)\bigr)^n= \E e^{tX_0},\qquad -\infty<t<\infty.\tag1.11
$$

Consequently,
$$
\P (X\ge  \E X+a)\le e^{-ta-t \E X} \E e^{tX}\le e^{-ta-t \E X_0} \E e^{tX_0},
\tag1.12
$$
and thus every Chernoff type bound for the binomial variable $X_0$ derived
from (1.1), applies also to $X$.

\proclaim{Theorem 2}
The bounds \rom{(1.2)--(1.10)} hold also when $X=\sum_1^n I_i$
where $I_i\sim \Be(p_i)$ are
independent indicator variables and $p= \E X/n$, $q=1-p$.
\qed
\endproclaim
\Remark{Remark 4}
We do not claim that the actual tail probability $\P (X_0\ge
 \E X+a)$ is larger than $\P (X\ge  \E X+a)$,
 and indeed this is in general false as
is shown by the example $n=2$, $p_1=1/5$, $p_2=3/5$,
where $\P (X\ge1)=17/25$ while
$X_0\sim \Bi (2,2/5)$ and $\P (X_0\ge1)=16/25$.
\endRemark

As mentioned above, the bounds (1.2) and (1.3) are asymptotically sharp for
the binomial distribution, but that is no longer generally true when  the
\0 variables have different distributions. In fact, a Taylor expansion
shows that the exponent in (1.2) or (1.3) is
$-\frac{a^2}{2np(1-p)}\bigl(1+o(1)\bigr)$
provided $a=o\bigl(np(1-p)\bigr)$, cf\. (1.4) and (1.7). In the binomial case, this
equals $-\frac{a^2}{2\sigma^2}\bigl(1+o(1)\bigr)$, with $\sigma^2=\Var X$,
which is what one
would expect from normal approximation heuristics; in general, however,
$\sigma^2=\Var X$ may be much smaller than $np(1-p)$, and it would be
advantageous to have better bounds with exponents
$-\frac{a^2}{2\sigma^2}\bigl(1+o(1)\bigr)$ for
moderately large $a$. 
This is achieved by {\it Bennett's inequality}, 
see Bennett (1962) and Hoeffding (1963),
which we state as (1.13) in the next
theorem; 
the simple consequence (1.14) is known as {\it Bernstein's inequality}, see
%Bernstein (1946) and 
Boucheron, Lugosi and Massart (2013).
Note that these inequalitites give  bounds depending on $a$ and $\sigma^2$
only, with an exponent of the expected order for $a=o(\sigma^2)$. 
We give a
proof in Section 3 using (1.1) as above, but doing a more careful
estimation of $ \E e^{tX}$ than (1.11).

\proclaim{Theorem 3} % {\rm (Bennett)}}
Let $X$ be a random variable and suppose that there exist independent \0
variables $I_i\sim \Be (p_i)$, $i=1,\dots,n$, such that $X \eqd \sum_1^n I_i$.
Let $\lambda= \E X$ and $\sigma^2=\Var X$. Then
\roster
\item If $a\ge0$, then
$$
\P (X\ge  \E X+a)
\le \exp \biggl(-\sigma^2\Bigl(\Bigl(1+\frac{a}{\sigma^2}\Bigr)\ln
\Bigl(1+\frac{a}{\sigma^2}\Bigr)-\frac{a}{\sigma^2}\Bigr)\biggr)\tag1.13
$$
\item
If $a\ge0$, then
$$
\P (X\ge  \E X+a)
\le \exp\biggl(-\frac{a^2}{2\sigma^2}\Big/\Bigl(1+\frac{a}{3\sigma^2}
\Bigr)\biggr)
\le \exp\biggl(-\frac{a^2}{2\sigma^2}\Bigl(1-\frac{a}{3\sigma^2}
\Bigr)\biggr).
\tag1.14
$$
\item If $a\ge c\sigma^2$, with $c>0$, then
$$
\P (X\ge  \E X+a)\le \exp\bigl(-((1+c^{-1})\ln (1+c)-1)a\bigr).\tag1.15
$$
\endroster
The same estimates hold for $\P (X\le  \E X-a)$, and thus $\P (|X- \E X|\ge a)$ may be
estimated by twice the right hand sides in  \rom{(1.13)--(1.15)}.
\roster
\item"{\rm(iv)}" Moreover, if\/ $a\ge0$ and $\sigma^2\ge\lambda/2$, then
$$
\P (X\le  \E X-a)\le \exp\Bigl(-\frac{a^2}{2\sigma^2}\Bigr).\tag1.16
$$
\endroster
\endproclaim
\Remark{Remark 5}
The estimates (1.13) and (1.14) are very similar to (1.5) and (1.6); the only
difference is that $np=\lambda$ is replaced by $\sigma^2<\lambda$.
It is easily seen that this always improves the bound in (1.5) and the first
bound in (1.6). On the other hand, the bounds for $\P (X\le  \E X-a)$
(except (1.16))
are
somewhat different from the corresponding bounds in Theorem 1, because of the
symmetry of the bounds in Theorem 3.
\endRemark
\Remark{Remark 6}
It is easily seen (by approximating a Poisson distribution) that the constant
in (iii) is best possible. In particular, it follows that
(1.16) cannot hold without restriction.
\endRemark
Much more precise estimates of the tail probabilities for sums of independent,
but not necessarily identically distributed, random variables were obtained by
Feller (1943) using
different methods (conjugated distributions as in
Cram\'er (1938) together with a Berry--Esseen estimate), and
it is interesting to compare our result with Feller's.
Feller's result (for our case, using $\lambda_n=1/\sigma$ in Feller (1943))
is, for  $0<a<\sigma^2/12$,
$$
\P (X\ge  \E X+a)=e^{-\frac{x^2}{2} Q(x)}
\Bigl(1-\Phi (x)+ \frac{\theta(x)}{\sigma}e^{-x^2/2}\Bigr),
\tag1.17
$$
where $x=a/\sigma$, $|\theta(x)|<9$, $\Phi$ is the normal distribution function
and
$$
Q(x)=\sum_{\nu=1}^\infty q_\nu x^\nu,\tag1.18
$$
where $q_\nu$ depends on the first $\nu+2$ moments of $X$ and
$$
|q_\nu|< \frac{1}{7}\Bigl(\frac{12}{\sigma}\Bigr)^\nu.\tag1.19
$$

If, say, $\sigma\le a\le \sigma^2/24$, this yields, using $1-\Phi (x)\le
(2\pi)^{-1/2} x^{-1}e^{-x^2/2}$ for $x>0$,
$$
\align
\P(X\ge  \E  X +a)&< \exp\Bigl(-\frac{x^2}{2}\bigl(1+Q(x)\bigr)\Bigr)
=\exp\biggl(-\frac{a^2}{2\sigma^2}\Bigl(1+Q\bigl(\frac{a}{\sigma}\bigr)\Bigr)\biggr)\\
&< \exp
\Bigl(-\frac{a^2}{2\sigma^2}\Bigl(1-\frac{24}{7} \frac{a}{\sigma^2}\Bigr)\Bigr).\tag1.20
 \endalign
$$
For $\P(X\le \E X-a)$ we have the same estimates if we replace $Q(x)$ by
$Q(-x)$ (and $\theta (x)$ by some $\theta'(x))$;
this follows by considering $n-X$.

The bound (1.20) is similar to the ones given in Theorem 3, in particular
(1.14). It is somewhat inferior to (1.14) since the constant in the second
order term in the exponent is worse, and the range of $a$ is restricted, but
for applications they are essentially equivalent.

Note also that Feller's result has other advantages. First, (1.17) is an
equality (although the exact value of $\theta(x)$ is unspecified), and it leads
also to a lower bound
similar to (1.20) and to asymptotic results.
In particular, simple asymptotic results follow when $a/\sigma^2\to0$ and thus
$Q(x)\to0$.
Secondly, Feller (1943) describes
how the coefficients $q_\nu$ may be explicitly expressed in terms of the
semi-invariants $\kappa_j$ of $X$ (and thus in terms of the moments); for
example (the sign seems to be wrong in Feller (1943), (2.18)--(2.19)),
$$
\align
q_1 &= -\frac{\kappa_3}{3\sigma^3},\tag1.21\\
q_2 &= -\frac{\kappa_4}{12\sigma^4}+\frac{1}{4\sigma^6}\kappa_3^2.\tag1.22
\endalign
$$
(Thus, $q_1=-\frac{1}{3}\gamma_1$ and
$q_2=-\frac{1}{12}\gamma_2+\frac{1}{4}\gamma_1^2$, where $\gamma_1$ and
$\gamma_2$ are the skewness and excess of $X$, respectively.)
For example, using (1.21) for $q_1$ and (1.19) for $q_\nu$, $\nu\ge2$,
we obtain
for $\sigma\le a\le \sigma^2/24$, instead of (1.20),
$$
\P (X\ge \E X+a)\le \exp
\biggl(-\frac{a^2}{2\sigma^2}\Bigl(1-\frac{288}{7}\frac{a^2}{\sigma^4}-\frac{a\kappa_3}
{3\sigma^4} \Bigr)\biggr),\tag1.23
$$
which yields an improvement in cases when $\kappa_3=\E (X-\E X)^3$
is known and either negative or not to large positive.
%\endproclaim
\subhead
Dependent summands
\endsubhead

Let us now consider the case of dependent \0 variables $I_i$. Of course any
bounded non-negative integer valued random variable $X$ can be written as a
sum of dependent \0 variables, so nothing can be said in general. We will
here consider only \0 variables that are {\it negatively related\/}
in the following sense, cf\. Barbour, Holst and Janson (1992).

(Note that large deviation bounds for a class of sums of {\it positively\/}
related indicators are given in Janson (1990) and Barbour, Holst and Janson
(1992), Theorem 2.S. In this case only the lower tail probabilities $\P (X\le
\E X-a)$ have nice upper bounds.)

\definition{Definition}
The indicator random variables $(I_i)_{i=1}^n$ (defined on the same
probability space) are negatively related if for each $j\le n$ there exist
further random variables $(J_{ij})_{i=1} ^n$, defined on the same probability
space (or an extension of it), such that the distribution of the random vector
$(J_{ij})_{i=1}^n$ equals the conditional distribution of $(I_i)_{i=1}^n$
given $I_j=1$, and, moreover, for every $i$  with $i\ne j$,
$J_{ij}^i \le I_i$.
\enddefinition
\example{Example 1} (Hypergeometric distribution.)
Let $m$, $n$ and $N$ be given positive integers with $\max(m,n)\le N$.
%Distribute $m$ balls into $N$ urns at random by drawing without replacement.
Given $N$ urns, labelled $1,\dots,N$, and $m$ balls, put the balls at random
into $m$ different urns (drawing without replacement), and let $X$ be the
total number of balls in urns $1,\dots,n$. Clearly $X=\sum_1^n I_i$, where
$I_i$ equals 1 if urn $i$ contains a ball. In this case it is easy to show
that the indicators $I_i$ are negatively related by explicitly construction
$J_i$, as follows. After randomly distributing the balls as above, and
recording $I_i$, we ensure that there is a ball in urn $j$ by ``cheating'': if
urn $j$ is empty we select one of the balls at random and move it to urn $j$.
Let $J_{ij}=1$ if urn $i$ now contains a ball. It is clear that $(J_{ij})$ has
the right distribution, and that $J_{ij}\le I_i$ for $i\ne j$.
\endexample
\example{Example 2} Distribute $m$ balls into $n$ urns, but this time put the
balls one by one at random, independently of the other choices of urn (drawing
with replacement). Let $X$ be the number of empty urns. Clearly $X=\sum_1^n
I_i$, where $I_i=1$ if urn $i$ is empty. These indicators are negativlely
related; this follow by a construction very similar to the one in Example 1,
removing all balls (if any) in urn $i$ and redistributing them (repeating if
necessary).
\endexample
Further examples of negatively related variables are given in Barbour, Holst
and Janson  (1992), where also some general results are established. In
particular, it is proven (a special case of Corollary 2.D.1) that the
variables
$I_i$, $i=1,\dots,n$, are negatively related if and only if $I_j$ and
$\phi(I_1,\dots,I_{j-1},I_{j+1},\dots,I_n)$ are negatively correlated for
every $j$ and every indicator function $\phi $ that is increasing in each
variable. (Pairwise negative correlation of the $I_i$ is not enough.) It
follows immediately that the variables $(1-I_i)_{i=1}^n$ are negatively
related if $(I_i)_{i=1}^n$ are. It follows also that variables
are negatively
related if they are
negatively associated in the sense of Joag-Dev and Proschan (1983).
\proclaim{Theorem 4}
Suppose that $X\eqd\sum_1^n I_i$, where $I_i\sim \Be (p_i)$ are negatively
related indicator variables. Let $\widetilde I_i$, $i=1,\dots,n$, be independent
indicator variables with $\widetilde I_i\sim \Be (p_i)$, and put $\widetilde
X=\sum_1^n \widetilde I_i$. Then, for every real $t$,
$$
 \E e^{tX}\le  \E e^{t\widetilde X}.
$$
Consequently, any Chernoff type bound for $\widetilde X$ applies also to $X$. In
particular, \rom{(1.2)--(1.10)} hold with $p=\E X/n$ and $q=1-p$; for example, with
$\lambda=\E X$,
$$
\align
\P  (X\ge \E X+a) & \le \exp \biggl(- \frac{a^2}{2\lambda (1+a/3\lambda)}\biggr)
\le \exp \biggl(-\frac{a^2}{2\lambda}\Bigl(1-\frac{a}{3\lambda}\Bigr)\biggr),\\
\P (X\le \E X-a) &\le \exp \biggl(-\frac{a^2}{2\lambda}\biggr).
\endalign
$$
\endproclaim

Of course, also the bounds in Theorem 3 (applied to $\widetilde X$) apply to $X$.
The problem is that we have to use $\sigma^2=\Var (\widetilde X)$ instead of $\Var
(X)$, which may be much smaller. In fact, the bounds in Theorem 3 are in
general false with $\sigma^2=\Var X$ in the dependent case; the following
theorem implies that it is impossible to have a general bound that is, say,
$\exp (- a^2/3\sigma^2)$ when $a= 4\sigma$.

\proclaim{Theorem 5}
Let $\alpha>0$, $0<c<1/e$ and $A<\infty$. There exists a random variable $X$
which is a finite sum of negatively related indicators such that
$\sigma^2=\Var X>A$ and, with $a=\alpha\sigma$,
$$
\P(X>\E X+a)>ce^{-a/\sigma}. \tag1.24
$$
\endproclaim

Nevertheless, there are cases where it is possible to do better. A striking
example is based on the result by Vatutin and Mikhailov (1982) that certain
random variables that occur in some occupancy problems, and have natural
representations as sums of negatively related indicators (with the same
expectation), also can be represented as sums of {\it independent\/}
indicators with different expectations. (The proof is algebraic, and based on
showing that the probability generating function has only real roots; there is
no (known) probabilistic interpretation of these indicators, which in general
have irrational expectations.) Their result includes the variables in Examples
1 and 2 (using in their notation $s_1=N-n$, $s_2=N-m$ for Example 1 and
$s_1=\dots =s_m=1$ for Example 2).

Consequently, the variables in Examples 1 and 2 actually satisfy the hypothesis
of Theorem 3 (although we do not know the $p_i$ explicitly). Hence we can
apply Theorem 3; note that the bounds in Theorem 3 involve only $\sigma^2$ and
possibly $\lambda$, and not the unknown $p_i$.
(In fact, this application was one of
the motivations for finding bounds of the form given in Theorem 3.)
\proclaim{Theorem 6}
Let $X$ be either hypergeometric as  in Example \rom1, or as in Example \rom2.
Then the conclusions of Theorem \rom3 hold, with $\lambda=\E X$ and
$\sigma^2=\Var X$.
\qed
\endproclaim

\example{Example 2, cont}
For the occupancy problem described above,
$$\gather
\lambda=\E X =n\Bigl(1-\frac1n\Bigr)^m,\\
\sigma^2=\Var X =n\Bigl(1-\frac1n\Bigr)^m + n(n-1)\Bigl(1-\frac2n\Bigr)^m
- n^2\Bigl(1-\frac1n\Bigr)^{2m};
\endgather
$$
estimates of the tail probabilities are obtained by using these values
in any of the formulas (1.2)--(1.16), letting $p=\E X/n$ and $q=1-n$.

For asymptotical results in the case $m/n\to r>0$, we easily find
$$\gather
\lambda=\E X \sim ne^{-r},\\
\sigma^2=\Var X \sim n \bigl(e^{-r}-(1+r)e^{-2r}\bigr) = ne^{-2r}(e^r-1-r).
\endgather
$$
The asymptotics for the tail probabilites in this case have been studied in
detail by Kamath, Motwani, Palem and Spirakis (1994).
\endexample

\Remark{Remark 7}
A comparison of Theorems 3 and 5 shows that not every random
variable that is a sum of negatively related indicators can be represented as
a sum of independent indicators; the Vatutin--Mikhailov result depends on some
further structure. The first example of such a variable was found by Andrew
Barbour (personal communication): Let $\P
(X=3)=4/13$, $\P(X=4)=5/13$, $\P(X=5)=4/13$. Then $X=\sum_1^5 I_i$, with $I_i$
indicators and the distribution of $(I_1,\dots,I_5)$ uniform given $X$; and
these $I_i$ are easily verified to be negatively related. On the other
hand, it is easily seen that $X$ is not the sum of any number of independent
indicators, since the probability generating function has non-real roots.
\endRemark

\subhead{Acknowledgements}
\endsubhead
This research has been inspired by discussions with
Andrew Barbour,
Carl-Gustav Esseen,
Rajeev Motwani,
Joel Spencer,
Andrew Thomason,
and possibly others.

\head
2. Proof of Theorem 1
\endhead

The first inequality in (1.4) is trivial for $a>nq$. For $0\le a\le nq$,
let $x=a/n \in [0,q]$. Then the bound (1.2) may be written
$$
\P(X\ge\E X+a)\le \exp\biggl(-np\Bigl(1+\frac{x}{p}\Bigr)\ln
\Bigl(1+\frac{x}{p}\Bigr)
-nq\Bigl(1-\frac{x}{q}\Bigr)\ln \Bigl(1-\frac{x}{q}\Bigr)\biggr).
\tag2.1
$$
Let, for $0\le x\le q$,
$$
f(x)= p\Bigl(1+\frac{x}{p}\Bigr)\ln \Bigl(1+\frac{x}{p}\Bigr)+ q\Bigl(1-\frac{x}{q}\Bigr)\ln
\Bigl(1-\frac{x}{q}\Bigr)-\frac{x^2}{2\bigl(pq+x(q-p)/3\bigr)}.
$$
Then $f(0)=f'(0)=0$, and an elementary calculation yields
$$\align
f''(x)&=\frac{1}{x+p}+\frac{1}{q-x}-\frac{p^2q^2}{\bigl(pq+x(q-p)/3\bigr)^3}\\
&=
\frac{\frac{1}{3} pq(q-p)^2x^2+\frac{1}{27} (q-p)^3x^3+p^2q^2x^2}
{(x+p)(q-x)\bigl(pq+x(q-p)/3\bigr)^3}\ge0
\endalign
$$
for $0\le x\le q$. Hence $f(x)\ge 0$ in this interval, and thus
$$
\P( X\ge\E X+a)\le \exp \biggl(-n \frac{x^2}{2\bigl(pq+x(q-p)/3\bigr)}\biggr).
$$
This proves the first inequality in (1.4). The second follows from $q-p\le 1$
and the third from
$$
\frac{1}{1+a/3npq} \ge 1-\frac{a}{3npq}.
$$

Inequality (1.5) follows directly from (2.1) and
$$
-nq\Bigl(1-\frac{x}{q}\Bigr)\ln \Bigl(1-\frac{x}{q}\Bigr)
=n(q-x)\ln\Bigl(1+\frac{x}{q-x}\Bigr) \le nx.
\tag2.2
$$

The inequalities (1.6) follow from (1.4) and
$$
-\frac{1}{npq+a/3} \le -\frac{1}{np(1+a/3np)}\le -\frac{1}{np}
\Bigl(1-\frac{a}{3np}\Bigr);
$$
alternatively, they follow easily from (1.5), cf\. (3.10).

The inequalities (1.7) follow from (1.3) by an argument similar to the one
given above for (1.4), or (simpler) by applying (1.4) to
$n-X\sim \Bi (n,q)$;
(1.8) follows from (1.3), using (2.2) with $x=-a/n\le0$;
(1.9) follows from (1.7) and, assuming (as we may) $a\le np$,
$npq-a(q-p)/3\le npq+ap/3\le npq+np^2/3\le np$. Finally, also (1.10) follows
from (1.7) since we now assume $a(q-p)/3\ge 0$.
\qed

\head{3. Proof of Theorem 3}
\endhead
We may assume that $X=\sum_1^n I_i$ where $I_i\sim \Be (p_i)$ are independent.
Note that
$$
\align
\lambda &= \E X=\sum_1^n p_i\\
\sigma^2 &= \Var X=\sum_1^n p_i(1-p_i)=\lambda-\sum_1^n p_i^2
\endalign
$$
and thus
$$
\sum_1^n p_i^2=\lambda-\sigma^2.
$$
We assume, to avoid trivialities, that at least one $p_i\ne 0,1$. Thus
$0<\lambda<n$ and $0<\sigma^2<\lambda$.

We begin with a real analysis lemma. It is an analogue of Jensen's inequality
but with a condition on the sign of the third derivative instead of the
second.

\proclaim{Lemma 1}
Suppose that $\mu$ is a finite positive measure on $[0,1]$, and define
$$
\align
m &= \mu([0,1]),\\
x_0 &= \int_0^1 x^2d\mu\bigg/\int_0^1 x\,d\mu,\\
\alpha_0 & =\biggl(\int_0^1 x\,d\mu\biggr)^2\bigg/\int_0^1x^2d\mu,\\
x_1 &= 1-\int_0^1 (1-x)^2 d\mu\bigg/\int_0^1 (1-x)d\mu,\\
\alpha_1 &= \biggl(\int_0^1(1-x)d\mu\biggr)^2 \bigg/\int_0^1(1-x)^2d\mu.
\endalign
$$
\rom(We here let $0/0=0$; this occurs in the degenerate cases where $\mu$ is
a point
mass at \rom0 or \rom1.\rom)
If $f$ is a three times continuously differentiable real function on $[0,1]$
with $f'''\ge0$, then
$$
(m-\alpha_0)f(0)+\alpha_0 f(x_0)\le \int_0^1 f\,d\mu\le
(m-\alpha_1)f(1)+\alpha_1f(x_1).\tag 3.1
$$
If instead $f'''\le0$ on $[0,1]$, then these inequalities are reversed.
\endproclaim

\demo{Proof}
We will show the left inequality of (3.1); the right inequality then follows by
symmetry, considering the function $\widetilde f(x)=-f(1-x)$, which satisfies
$\widetilde f'''(x)=f'''(1-x)\ge0$, and the similarly reflected measure $\widetilde
\mu(A)=\mu(\{1-x:x\in A\})$.
Similarly, the statement for $f'''\le 0$ follows by considering $-f$.

Let $\nu$ be the measure $(m-\alpha_0)\delta_0+\alpha_0\delta_{x_0}$; thus the
sought inequality is $\int f\,d\nu\le \int f\,d\mu$, while the choice of
$\alpha_0$ and $x_0$ yields
$\int  1\, d\nu=m=\int 1\,d\mu$, $\int x\,d\nu=\alpha_0x_0=\int x\,d\mu$, and
$\int
x^2 d\mu=\alpha_0x_0^2=\int x^2d\mu$. (In fact, it is easily seen that $\nu$
is the unique measure concentrated on a two-point set $\{0,x\}$ for some $x\in
[0,1]$, such that $\int x^k d\nu=\int x^k d\mu$ for $k=0,1,2$.)

We now use Taylor's formula
$$
f(x)=f(0)+f'(0)x+\tfrac{1}{2}f''(0)x^2+\tfrac{1}{2}\int_0^x (x-t)^2f'''(t)\,dt
\tag 3.2
$$
and integrate against the signed measure $\mu-\nu$. Since as we just have
shown,
$$
\int_0^1 x^k d(\mu-\nu)=0,\qquad k=0,1,2,\tag 3.3
$$
we obtain from (3.2) and Fubini's theorem,
$$
\align
\int_0^1 f(x)\,d(\mu-\nu) &= \int_0^1\int_0^x \tfrac{1}{2}(x-t)^2
f'''(t)\,dt\,d(\mu-\nu)(x)\\
&= \int_0^1 \int_t^1 \tfrac{1}{2}(x-t)^2 d(\mu-\nu)(x)f'''(t)\,dt\\
&= \int_0^1 \varphi (t)f'''(t)dt,
\tag {3.4}
\endalign
$$
where
$$
\varphi(t)=\tfrac{1}{2}\int_t^1 (x-t)^2 d(\mu-\nu)(x).
$$
We claim that $\varphi(t)\ge0$ on [0,1]; this implies $\int_0^1
f(x)d(\mu-\nu)\ge0$ by (3.4) which is the required result.
Note that $\varphi(1)=0$ and $\varphi(0)=\frac{1}{2}\int_0^1 x^2d(\mu-\nu)=0$.
Moreover, again by (3.3), $\int_0^1 (x-t)^2d(\mu-\nu)(x)=0$ and thus
$$
\varphi(t)=-\frac{1}{2}\int_0^t(x-t)^2d(\mu-\nu)(x).
$$
 Using Fubini again, and letting $F(x)=\mu([0,x])-\nu([0,x])$,
 $$
 \varphi(t)=-\iiint_{0\le x\le y\le z\le t}   dy\,dz
 \,d(\mu-\nu)(x)=-\int_0^t \int _0^zF(y)\,dy\,dz,
 $$
 and thus $\varphi$ is continuously differentiable with
 $$
 \varphi'(z)=-\int_0^z F(y)\,dy.\tag 3.5
 $$
 In particular $\varphi'(0)=0$ and, using Fubini a last time,
 $$
 \varphi'(1)=-\int_0^1 F(y)\,dy=-\int_0^1\int_0^y d(\mu-\nu)(x)\,dy=-\int_0^1
 (1-x)\,d(\mu-\nu)(x)=0.
 $$
 On the interval $[0,x_0)$, $\nu([0,x])$ is constant $m-\alpha_0$, and thus
 $F(x)$ is increasing; hence there exists $x_2\in [0,x_0]$ such that
 $F(x)\le0$ on $(0,x_2)$ and $F(x)\ge0$ on $(x_2,x_0)$. It follows by (3.5)
 that $\varphi'$ is increasing and thus $\varphi$ is convex on $[0,x_2]$,
 while $\varphi$ is concave on $[x_2,x_0]$.
 Since $\varphi(0)=\varphi'(0)=0$, this implies that $\varphi\ge0$ on
 $[0,x_2]$.
Similarly,
 on the interval $[x_0,1]$, we have $\nu([0,x])=m$ and thus
 $$
 F(x)=\mu([0,x])-\nu([0,x])=\mu([0,x])-m\le 0,
 $$
which implies that $\varphi$ is convex on $[x_0,1]$.
Moreover, $\varphi(1)=\varphi'(1)=0$ and thus $\varphi\ge0$ on
 $[x_0,1]$. Finally, on the interval $[x_2,x_0]$, $\varphi$ is concave so it
 attains it minimum at one of the endpoints, but we have already shown
 $\varphi(x_2),\varphi(x_0)\ge0$ and thus $\varphi\ge0$ also on $[x_2,x_0]$,
 which completes the proof.
 \qed
\enddemo

We apply this lemma to estimate the moment generating function of $X$.

\proclaim{Lemma 2}
Let $X$ be as above. If\/ $0\le t\le 1$, then
$$
\E (1-t)^X\le \biggl(1-t\Bigl(1-\frac{\sigma^2}{\lambda}\Bigr)\biggr)^{\lambda^2/(\lambda-\sigma^2)}
$$
or
$$
\ln \E(1-t)^X\le \frac{\lambda^2}{\lambda-\sigma^2}\ln
\biggl(1-t\Bigl(1-\frac{\sigma^2}{\lambda}\Bigr)\biggr).
\tag 3.6
$$
\endproclaim

\demo{Proof}
Since the $I_i$ are independent, and $\E(1-t)^{I_i}=1-p_i+p_i(1-t)=1-p_i t$,
$$
\E(1-t)^X =\E\prod_{i=1}^n(1-t)^{I_i}=\prod_{i=1}^n \E(1-t)^{I_i}=\prod_{i=1}^n(1-p_it),
$$
and thus
$$
\ln \E(1-t)^X=\sum_{i=1}^n\ln (1-p_it)=\int_0^1\ln (1-tx)\,d\mu(x),
$$
where $\mu$ is the measure $\sum_1^n \delta_{p_i}$ consisting of $n$
point masses at the (possibly coinciding) points $p_i$. Note that
$$
\int_0^1x\,d\mu=\sum_{i=1}^np_i=\lambda
$$
and
$$
\int_0^1x^2d\mu=\sum_{i=1}^np_i^2=\lambda-\sigma^2.
$$
We may assume that $t<1$ (the case $t=1$ follows then by continuity); then the
function $f(x)=\ln(1-tx)$ is infinitely differentiable on $[0,1]$ with
$f'''(x)=-2t^3/(1-tx)^3\le0$. Hence Lemma 1 yields
$$
\ln \E(1-t)^X=\int_0^1\ln
(1-tx)\,d\mu(x)\le(m-\alpha_0)f(0)+\alpha_0f(x_0)=\alpha_0\ln(1-tx_0),
$$
where $m=n$, $\alpha_0=\lambda^2/(\lambda-\sigma^2)$ and
$x_0=(\lambda-\sigma^2)/\lambda$, which is the required estimate.
\qed
\enddemo
\Remark{Remark 8}
For $t\ge0$, a similar argument yields
$$
\align
\ln\E(1+t)^X & \le(n-\alpha_1)\ln(1+t)+\alpha_1\ln(1+tx_1)\\
&=\frac{n\lambda-\lambda^2-n\sigma^2}
{n-\lambda-\sigma^2}\ln(1+t)+\frac{(n-\lambda)^2}{n-\lambda-\sigma^2}\ln\Bigl(1+t
\frac{\sigma^2}{n-\lambda}\Bigr).\tag 3.7
\endalign
$$
This inequality could be used instead of (3.6) below, giving the same results.
We prefer to use (3.6), which does not involve $n$ explicitly.
\endRemark
\Remark{Remark 9}
Estimates of $\E e^{sX}$ are, of course, obtained by substituting $t=1-e^s$ in
(3.6) for $s\le0$ and $t=e^s-1$ in (3.7) for $s\ge0$.
\endRemark

We can now obtain our basic estimate.
\proclaim{Lemma 3}
Let $X$ be as above.
If\/ $0\le a\le\lambda$, then
$$
\ln \P (X\le
\lambda-a)\le-\frac{\lambda}{\lambda-\sigma^2}(a+\sigma^2-a\sigma^2/\lambda)
\ln\Bigl(1+\frac{a}{\sigma^2}-\frac{a}{\lambda}\Bigr)
-(\lambda-a)\ln\Bigl(1-\frac{a}{\lambda}\Bigr).
\tag 3.8
$$
\rom(When $a=\lambda$, we define $(\lambda-a)\ln (1-a/\lambda)=0$.\rom)
\endproclaim
\demo{Proof}
For any $t$ with $0\le t\le1$,
$$
\E(1-t)^X\ge(1-t)^{\lambda-a}\P(X\le \lambda-a),
$$
and thus, using Lemma 2,
$$
\align
\ln \P (X\le\lambda-a)&\le\ln \E(1-t)^X-(\lambda-a)\ln(1-t)\\
&\le\frac{\lambda^2}{\lambda-\sigma^2}\ln\biggl(1-t\Bigl(1-\frac{\sigma^2}
{\lambda}\Bigr)\biggr)-(\lambda-a)\ln(1-t).
\endalign
$$
Choosing
$$
t=\frac{a}{a+\sigma^2-a\sigma^2/\lambda}
$$
(which minimizes the right  hand side), this yields
$$
\align
\ln\P (X\le \lambda-a)
&\le \frac{\lambda^2}{\lambda-\sigma^2}\ln
  \frac{\sigma^2}{a+\sigma^2-a\sigma^2/\lambda}
 -(\lambda-a)\ln  \frac{\sigma^2(1-a/\lambda)}{a+\sigma^2-a\sigma^2/\lambda}\\
&= \Bigl(\frac{\lambda^2}{\lambda-\sigma^2}-\lambda+a\Bigr)\ln
\frac{\sigma^2}{a+\sigma^2-a\sigma^2/\lambda}-(\lambda-a)\ln\Bigl(1-\frac{a}
{\lambda}\Bigr)\\
&=-\frac{\lambda\sigma^2+\lambda a-a\sigma^2}{\lambda-\sigma^2}\ln
\frac{\sigma^2+a-a\sigma^2/\lambda}{\sigma^2}
 -(\lambda-a)\ln\Bigl(1-\frac{a}{\lambda}\Bigr),
\endalign
$$
which yields the sought result.
\qed
\enddemo

While the estimate in Lemma 3 may be useful for numerical evaluation in
applications, it is too complicated to be of much other direct use.
Hence we will use it to derive the simpler (but slightly weaker) estimates in
Theorem 3.

For notational convenience, let $\sigma^2=x\lambda$ and
$a=y\sigma^2=xy\lambda$, where $0<x<1$ and $0\le y\le 1/x$. Then (3.8) may be
written
$$
\ln \P (X\le \lambda-a)\le -\sigma^2 g(x,y),\tag 3.9
$$
where
$$
g(x,y)=\frac{1}{1-x}(1+y-xy)\ln (1+y-xy)+\frac{1}{x}(1-xy)\ln (1-xy).
$$

\proclaim{Lemma 4}
$g(x,y)$ is an increasing function of $x$ in the region
$$
U=\{(x,y):0<x<1,\;0\le y\le 1/x\}.
$$
\endproclaim

\demo{Proof}
We want to show that $\partial g/\partial x\ge0$ in the region
$\{(x,y):0<x<1,\, 0\le y< 1/x\}$; note that $g$ is well-defined and infinitely
differentiable in the larger region $V=\{(x,y):0<x<1,\; -1/(1-x)<y<1/x\}$, and
continuous on $\bar V\cap \{(x,y):0<x<1\}\supset U$.

Instead of estimating $\partial g/\partial x$ directly, we first compute
$$
\frac{\partial g}{\partial y}=\ln (1+y-xy)-\ln (1-xy)
$$
and
$$
\frac{\partial ^2g}{\partial x\partial y}=\frac{\partial^2g}{\partial
y\partial x}=\frac{-y}{1+y-xy}-\frac{-y}{1-xy}=\frac{y^2}{(1+y-xy)(1-xy)}\ge0.
$$
Hence $\partial g/\partial x$ is an increasing function of $y$.
Moreover, taking $y=0$, we find $g(x,0)=0$ and thus
$$
\frac{\partial g}{\partial x}(x,0)=0,\qquad0<x<1.
$$
Hence $\frac{\partial g}{\partial x}(x,y)\ge0$ for all $(x,y)\in V$ with
$y\ge0$.
\qed
\enddemo
\demo{Proof of Theorem 3}
If $x\searrow 0$ and $y\ge0$ is fixed, then, as is easily seen,
$$
g(x,y)\to(1+y)\ln(1+y)-y.
$$
Lemma 4 thus yields, for $0<x<1$ and $0\le y\le 1/x$,
$$
g(x,y)\ge(1+y)\ln(1+y)-y,
$$
and thus by (3.9), for $0\le a\le\lambda$,
$$
\ln \P (X\le \lambda-a)\le -\sigma^2\bigl((1+y)\ln(1+y)-y\bigr),
$$
with $y=a/\sigma^2$, which is the analogue of (1.13) for $\P (X\le \E X-a)$;
note that this estimate trivially holds for $a>\lambda$.
In order to obtain (1.13), we consider the variable
$X^*=n-X\eqd \sum_1^n(1-I_i)$, and observe that
$\P(X\ge\E X+a)=\P(X^*\le \E X^*-a)$ and
$\Var (X^*)=\Var X$. The estimates (1.14) and (1.15) and their analogues for
$\P(X\le \E X-a)$ now follow from the elementary estimates, defining
$h(y)=(1+y)\ln(1+y)-y$,
$$
h(y)\ge\frac{y^2}{2(1+y/3)}\ge\frac{y^2}{2}(1-y/3),\qquad y\ge0,\tag 3.10
$$
and
$$
h(y)\ge yh(c)/c,\qquad y\ge c.\tag 3.11
$$
The estimate (3.10) may be verified by observing that $h(y) -y^2/2(1+y/3)$
vanishes together with its first derivative at 0, while the second derivative
equals $(9y^2+y^3)/(1+y)(3+y)^3\ge0$. Similarly, (3.11) follows by the
convexity of $h$. We omit the details.

Finally, if $\sigma^2 \ge \lambda/2$, then $x\ge 1/2$ and thus by Lemma 4
(assuming as we may that $a\le\lambda$),
$$
g(x,y) \ge g(\tfrac12,y) =(2+y)\ln(1+\tfrac y2) + (2-y)\ln(1-\tfrac y2)
\ge \tfrac12 y^2,
$$
which together with (3.9) yields (1.16).
\qed
\enddemo

\head
4. Proof of Theorem 4.
\endhead

Let $Y=\sum_2^n I_i$ and $Z=\sum_2^n J_{i1}$, where $(J_{i1})_{i=1}^n$ are as
in the definition of negatively related variables. Then $X=I_1+Y$ and 
$(Y\mid 
I_1=1){\mathop{\,=\,}\limits^d}Z$. Moreover, since $J_{i1} \le I_i$ for $i
\ge 2$, we have $Z \le Y$. Consequently, for any real $t$,
$$
\align
\E e^{tX} - \E e^{tY} & = \E(e^{tI_1}-1)e^{tY}=\E(e^t-1)I_1e^{tY} \\
& =(e^t-1)p_1\E(e^{tY}\mid I_1=1) \\
& =(e^t-1) p_1 \E e^{tZ} \\
& \le (e^t-1) p_1 \E e^{tY} \\
& = \E(e^{tI_1}-1) \E e^{tY}
\endalign
$$
and thus
$$
\E e^{tX} \le \E e^{tI_1} \E e^{tY}.
$$
Induction yields
$$
\E e^{tX} \le \prod_1^n e^{tI_i}=\E e^{t\widetilde X}.
\eqno{\square}
$$

\head
5. Proof of Theorem 5.
\endhead
Given $n$, $p$, $k$, with $1 \le k \le n$ and $0 < p < 1$, let
$U \sim \Bi(n,p)$
and let $X_{npk}$ be the random variable $U$ conditioned on $U \ge k$.
Since $U = \sum_1^n I_i$, with $I_i \sim \Be(p)$ independent,
$X_{npk}=\sum_1^n(I_i\mid \sum_1^n I_i \ge k)$, and it follows from Barbour,
Holst and Janson (1992) Proposition 2.2.10 and Theorem 2.I that $X_{npk}$ is
a sum of negatively related indicators. We claim that, for any $\alpha$, $c$,
$A$ as
in Theorem 5, some variable $X_{npk}$ satisfies (1.24).

Suppose not. Then for
each $n$, $p$, $k$ either $\Var(X_{npk}) \le A$ or
$$
\P (X_{npk} - \E X_{npk} > \alpha \sigma) \le c e^{-\alpha}. \tag{5.1}
$$
Fix $p\in(0,1)$, let $q=1-p$, choose $\varepsilon$ with $0 < \varepsilon < q$, take
$k=\lfloor n(p+\varepsilon)\rfloor+1$, and let $n \to \infty$. 
Then, for fixed $i \ge 0$, with
$r=p({q-\varepsilon})/(p+\varepsilon)q<1$,
$$
\frac{\P(X_{npk}=k+i+1)}{\P(X_{npk}=k+i)}=\frac{p}{q} \cdot
\frac{n-k-i}{k+i+1}\to
\frac{p}{q} \cdot \frac{q-\varepsilon}{p+\varepsilon}=r
\qquad \text{as $n\to \infty$}, \tag{5.2}
$$
and
$$
\P(X_{npk}=k+i) \le r^i \P(X_{npk}=k) \le r^i. \tag{5.3}
$$
It follows that $X_{npk}-k$ converges in distribution to a random variable
$Y_r$ with geometric distribution $\Ge(1-r)$: $\P(Y_r=i)=(1-r)r^i$, $i \ge 0$.
Moreover, by (5.3), every moment $\E(X_{npk}-k)^m$ stays bounded, which
implies that the moments converge to the corresponding moments of $Y_r$.

The
variance of $Y_r$ equals $r/(1-r)^2$. If $\Var(Y_r) > A$, then also
$\Var(X_{npk})> A$ for large $n$, so by our assumption (5.1) holds and taking
the limit as $n \to \infty$ we obtain
$$
\P(Y_r-\E Y_r > \alpha \sqrt{\Var Y_r}) \le c e^{-\alpha}. \tag{5.4}
$$
Now, let $\varepsilon \to 0$ (keeping $p$ fixed). Then $r \to 1$ and
$\Var(Y_r)=r/(1-r)^2 \to \infty$, so (5.4) holds when $r$ is close to 1.
Moreover, it is easily seen that as $r \to 1, (1-r)Y_r$ converges in
distribution to an exponential variable $Z \sim \Exp(1)$, again with
convergence of all moments. Consequently we may take the limit again and
obtain from (5.4)
$$
\P(Z-\E Z > \alpha \sqrt{\Var Z}) \le c e^{-\alpha}. \tag{5.5}
$$
But $\E Z=\Var Z = 1$, so the left hand side of (5.5) equals $\P(Z >
1+\alpha)=e^{-1-\alpha} > c e^{-\alpha}$, and we have obtained a contradiction.
\qed

\enddocument
\end